\documentclass[10pt]{article}
\usepackage[utf8]{inputenc}
\usepackage[T1]{fontenc}
\usepackage{amsmath}
\usepackage{amsfonts}
\usepackage{amssymb}
\usepackage[version=4]{mhchem}
\usepackage{stmaryrd}
\usepackage{hyperref}
\hypersetup{colorlinks=true, linkcolor=blue, filecolor=magenta, urlcolor=cyan,}
\usepackage{graphicx}
\usepackage[export]{adjustbox}
\graphicspath{ {./images/} }

\title{\Large A Short Note on the Infinity Product Tan (z) Function }

\author{Carlos A. Pérez Aparicio.${ }^{1}$}
\date{}

\begin{document}
\maketitle

${ }^{1}$ COITIRM, Independent Researcher.

Email: \href{mailto:cperapa@gmail.com}{cperapa@gmail.com}

Abstract:We present a novel derivation of the infinite product of the tangent function, Tan(z), expressed in terms of trigonometric expressions including Euler’s Sinc function and Viète’s formula, along with their generalizations.


\section*{1. Introduction}
Infinite product (1) was previously published by Viète in 1593, based on geometric considerations. A detailed compendium, including analysis and evaluations of infinite products of elementary functions in terms of trigonometric functions, can be found in the classical literature

\begin{equation*}
\prod_{k=0}^{\infty} \cos \left(\pi 2^{-k-2}\right)=\frac{2}{\pi} \tag{1}
\end{equation*}

Viète's formula is the following infinite product of nested radicals

\begin{equation*}
\frac{1}{8} \sqrt{\frac{1}{2}(\sqrt{2}+2)(\sqrt{\sqrt{2}+2}+2)(\sqrt{\sqrt{\sqrt{2}+2}+2}+2)...}=\frac{2}{\pi} \tag{2}
\end{equation*}

Viète's formula may be obtained as a special case of a formula for the sinc function that has often been attributed to Leonhard Euler[1], more than a century later:

\begin{equation*}
\prod_{k=0}^{\infty} \cos \left(2^{-k-1} z\right)=\frac{\sin (z)}{z} \tag{3}
\end{equation*}
\noindent

expressing each term of the product on the left as a function of earlier terms using the half-angle formula:

\begin{equation*}
\cos \left(\frac{z}{2}\right)=\sqrt{\frac{1}{2}(\cos (z)+1)} \tag{4}
\end{equation*}
\noindent

gives Viète's formula.

This formula has been generalization obtained through use of both Chebyshev polynomial and Fourier transformation (see Nishimura [3] and Kent E. Morrison [2]) for integer $q \geq 2$. As

\begin{equation*}
\prod_{j=0}^{\infty} \frac{\sum_{n=1}^{q} \cos \left(\frac{(2 n-1) z}{(2 q)^{j+1}}\right)}{q}=\frac{\sin (z)}{z} \tag{5}
\end{equation*}
\noindent

This can be written in more compact way

\begin{equation*}
\prod_{k=0}^{\infty} \frac{\sin \left(z q^{-k}\right) \csc \left(z q^{-k-1}\right)}{q}=\operatorname{sinc}(z) \tag{6}
\end{equation*}

\section*{2. The Tan Function}
We introduce two new generalizations of Euler's infinite product of the tangent function in the style of sinc(z) type,  $\mathrm{q} \geq 2$.
\begin{align*}
&\prod_{k=0}^{\infty} 2^{-\left((q-1) q^k\right)} q^{-q^k (k (q-1)+q)} z^{(q-1) q^k} \\
&\quad \times \tan^{-q^{k+1}}\left(\frac{1}{2} z q^{-k-1}\right) \tan^{q^k}\left(\frac{z q^{-k}}{2}\right) \\
&= \frac{2 \tan \left(\frac{z}{2}\right)}{z}
\end{align*}

\noindent

And the equivalent convergent product

\begin{equation*}
\prod_{k=0}^{\infty} \frac{\tan \left(\frac{z q^{-k}}{2}\right) \cot \left(\frac{1}{2} z q^{-k-1}\right)}{q}=\frac{2 \tan \left(\frac{z}{2}\right)}{z} \tag{8}
\end{equation*}

\section*{3. Generalization.}

It is possible to generalize equations (7) and (8) so that   $n$ and $m$ are integers $ \{n, m\} \geq 0$ as follows:

\begin{align*}
&\prod_{k=m}^{n-1} 2^{-\left((q-1) q^k\right)} q^{-q^k (k (q-1)+q)} z^{(q-1) q^k} \\
&\quad \times \tan^{-q^{k+1}}\left(\frac{1}{2} z q^{-k-1}\right) \tan^{q^k}\left(\frac{z q^{-k}}{2}\right) \\
&= 2^{q^m-q^n} q^{m q^m-n q^n} z^{q^n-q^m} \\
&\quad \times \tan^{q^m}\left(\frac{z q^{-m}}{2}\right) \tan^{-q^n}\left(\frac{z q^{-n}}{2}\right)\tag{9}
\end{align*}

And

\begin{align*}
& \prod_{k=m}^{n-1} \frac{\tan \left(\frac{z q^{-k}}{2}\right) \cot \left(\frac{1}{2} z q^{-k-1}\right)}{q}  \tag{10}\\
& =q^{m-n} \tan \left(\frac{z q^{-m}}{2}\right) \cot \left(\frac{z q^{-n}}{2}\right)
\end{align*}

For brevity, we employ mathematical induction to prove the case of Formula (10) when $q=2$ and $m=0$

\textbf{Proof}: starting with

\begin{equation*}
\prod_{k=0}^{n-1} \frac{1}{2} \tan \left(2^{-k-1} z\right) \cot \left(2^{-k-2} z\right)=2^{-n} \tan \left(\frac{z}{2}\right) \cot \left(2^{-n-1} z\right) \tag{11}
\end{equation*}
Step 1: Base Case $n=1$ We need to show that the formula holds for $n=1$ : the left-hand side LHS becomes:

\begin{equation*}
\frac{1}{2} \tan \left(2^{-1-1} z\right) \cot \left(2^{-1-2} z\right)=\frac{1}{2} \tan \left(\frac{z}{2}\right) \cot \left(\frac{z}{4}\right) \tag{12}
\end{equation*}
\noindent

The right-hand side RHS becomes:

\begin{equation*}
2^{-1} \tan \left(\frac{z}{2}\right) \cot \left(2^{-1-1} z\right)=\frac{1}{2} \tan \left(\frac{z}{2}\right) \cot \left(\frac{z}{4}\right) \tag{13}
\end{equation*}

The LHS and RHS are equal for the base case.

Step 2: Iinductive hypothesis assume that the formula is true for some arbitrary positive integer $n$:

\begin{equation*}
\prod_{k=0}^{n-1} \frac{1}{2} \tan \left(2^{-k-1} z\right) \cot \left(2^{-k-2} z\right)=2^{-k} \tan \left(\frac{z}{2}\right) \tag{14}
\end{equation*}
\noindent

Step 3: Inductive step we need to prove that the formula holds for $n+1$ :

Consider the LHS for $\mathrm{n}+1$ :

\begin{equation*}
\prod_{k=0}^{n} \frac{1}{2} \tan \left(2^{-k-1} z\right) \cot \left(2^{-k-2} z\right) \tag{15}
\end{equation*}
\noindent

Multiply both sides by $\left(\frac{1}{2} \tan \left(2^{-n-1} z\right) \cot \left(2^{-n-2} z\right)\right)$ :

\begin{equation*}
\left(\prod_{k=0}^{n-1} \frac{1}{2} \tan \left(2^{-k-1} z\right) \cot \left(2^{-k-2} z\right)\right) \cdot \frac{1}{2} \tan \left(2^{-\mathrm{n}-1} z\right) \cot \left(2^{-\mathrm{n}-2} z\right) \tag{16}
\end{equation*}
\noindent

Using the inductive hypothesis, the RHS simplifies to:

\begin{equation*}
2^{-k} \tan \left(\frac{z}{2}\right) \cot \left(2^{-\mathrm{n}-1} z\right) \cdot \frac{1}{2} \tan \left(2^{-\mathrm{n}-1} z\right) \cot \left(2^{-\mathrm{n}-2} z\right) \tag{17}
\end{equation*}
\noindent
Simplify the RHS:

\begin{equation*}
2^{-\mathrm{n}-1} \tan \left(\frac{z}{2}\right) \cot \left(2^{-\mathrm{n}-1} z\right) \tag{18}
\end{equation*}

The RHS matches the LHS for $n+1$.

Therefore, by mathematical induction, the given formula holds for all positive integers (n)

Example 1. Infinitely nested radicals.

Using several values for $n$ Integers $\frac{2 n \tan \left(\frac{\pi}{2 n}\right)}{\pi}$ in formula (10) get the following table:
\begin{center}
\begin{tabular}{|c|c|}
\hline
Product $\infty$ & Result \\
\hline
$\frac{1}{4} \sqrt{\frac{2+\sqrt{2+\sqrt{2}}}{2-\sqrt{2+\sqrt{2}}}} \cdots$ & $\frac{4}{\pi}$ \\
\hline
$\frac{(2-\sqrt{3})(2+\sqrt{3})\left(\frac{1}{4} \sqrt{3(2-\sqrt{2})}(1+\sqrt{2})+\frac{1}{4}(-1+\sqrt{2}) \sqrt{2+\sqrt{2}}\right)}{4 \sqrt{3}\left(\frac{1}{4} \sqrt{2-\sqrt{2}}(1+\sqrt{2})-\frac{1}{4}(-1+\sqrt{2}) \sqrt{3(2+\sqrt{2})}\right)} \cdots$ & $\frac{2 \sqrt{3}}{\pi}$ \\
\hline
$\frac{1}{4} \sqrt{\frac{(2-\sqrt{2})(2+\sqrt{2+\sqrt{2+\sqrt{2}})}}{(2+\sqrt{2})(2-\sqrt{2+\sqrt{2+\sqrt{2}})}}} \cdots$ & $\frac{8 \sqrt{\frac{2-\sqrt{2}}{2+\sqrt{2}}}}{\pi}$ \\
\hline
$\frac{\sqrt{1-\frac{2}{\sqrt{5}}}\left(\frac{1}{8}(-1+\sqrt{2}) \sqrt{2+\sqrt{2}}(-1+\sqrt{5})+\frac{1}{4}(1+\sqrt{2}) \sqrt{\frac{1}{2}(2-\sqrt{2})(5+\sqrt{5})}\right)}{4\left(-\frac{1}{8} \sqrt{2-\sqrt{2}}(1+\sqrt{2})(-1+\sqrt{5})+\frac{1}{4}(-1+\sqrt{2}) \sqrt{\frac{1}{2}(2+\sqrt{2})(5+\sqrt{5})}\right)} \cdots$ & $\frac{10 \sqrt{1-\frac{2}{\sqrt{5}}}}{\pi}$ \\
\hline

 & $\frac{12(2-\sqrt{3})}{\pi}$ \\
\hline
$\frac{1}{4} \sqrt{\frac{(2-\sqrt{2+\sqrt{2}})(2+\sqrt{2+\sqrt{2+\sqrt{2+\sqrt{2}}})}}{(2+\sqrt{2+\sqrt{2}})(2-\sqrt{2+\sqrt{2+\sqrt{2+\sqrt{2}}})}}} \cdots$ & $\frac{16 \sqrt{\frac{2-\sqrt{2+\sqrt{2}}}{2+\sqrt{2+\sqrt{2}}}}}{\pi}$ \\
\hline
\end{tabular}
\end{center}
\begin{center}
Table 1. 
\end{center}

\section*{References}
\begin{enumerate}
    \item Euler, Leonhard (1738). "De variis modis circuli quadraturam numeris proxime exprimendi"
    \item Kent E. Morrison (1995) Cosine Products, Fourier Transforms, and Random Sums, \emph{The American Mathematical Monthly}, 102:8, 716-724, DOI: 10.1080/00029890.1995.12004647.
    \item R. Nishimura. A generalization of Viète's infinite product and new mean iterations. \emph{The Australian Journal of Mathematical Analysis and Applications}, 13(1), 2016.
\end{enumerate}

\end{document}